\documentclass[12pt]{article}
\usepackage{latexsym,amssymb,amsmath}
\usepackage{a4wide}
\newtheorem{thm}{Theorem}
\newtheorem{lem}{Lemma}
\newtheorem{cor}{Corollary}
\newtheorem{prop}{Proposition}
\newcommand{\ebox}{\hfill $\Box$\vspace{3mm}}
\newcommand{\pr}{{\bf Proof.}\ }
\newcommand{\bt}{\begin{thm}}
\newcommand{\et}{\end{thm}}
\newcommand{\bl}{\begin{lem}}
\newcommand{\el}{\end{lem}}
\newcommand{\bp}{\begin{prop}}
\newcommand{\ep}{\end{prop}}
\newcommand{\bc}{\begin{cor}}
\newcommand{\ec}{\end{cor}}
\newcommand{\be}{\begin{eqnarray}}
\newcommand{\ee}{\end{eqnarray}}
\newcommand{\bi}{\begin{itemize}}
\newcommand{\ei}{\end{itemize}}
\newcommand{\beq}{\begin{equation}}
\newcommand{\eeq}{\end{equation}}
\newcommand{\noi}{\noindent}

\begin{document}

\title{\vspace{-0.9cm} Inverting non-invertible trees}


\author{}
\date{}
\maketitle

\begin{center}
\vspace{-1.7cm}

{\large So\v{n}a Pavl\'ikov\'a} \\
\vspace{2mm} {\small Slovak University of Technology, Slovakia}

\vskip 4mm

{\large Jozef \v Sir\'a\v n} \\
\vspace{2mm} {\small Open University, U.K., and Slovak University of Technology, Slovakia}

\end{center}
\vskip 3mm

\begin{abstract}
If a graph has a non-singular adjacency matrix, then one may use the inverse matrix to define a (labeled) graph that may be considered to be the inverse graph to the original one. It has been known that an adjacency matrix of a tree is non-singular if and only if the tree has a unique perfect matching; in this case the determinant of the matrix turns out to be $\pm 1$ and the inverse of the tree was shown to be `switching-equivalent' to a simple graph [C. Godsil, Inverses of Trees, Combinatorica 5 (1985), 33--39]. Using generalized inverses of symmetric matrices (that coincide with Moore-Penrose, Drazin, and group inverses in the symmetric case) we prove a formula for determining a `generalized inverse' of a tree.

\vskip 3mm

\noi {\em Keywords:} Graph Inverse; Tree; Adjacency matrix; Maximum matching.
\vskip 1mm

\noi {\em AMS Subject Classification:} 05C50, 15A09, 05C05, 05C70.

\end{abstract}

\section{Introduction}\label{intro}

We will consider finite, undirected graphs with no multiple edges but we allow every vertex to carry at most one loop (an edge whose both ends are the same vertex); such an object will simply be referred to as a {\em graph} in this paper. We will further assume that each edge $e$ of a graph $G$ carries a non-zero real {\em label} $\alpha(e)$, and the pair $(G,\alpha)$ will be called an {\em edge-labeled graph}. Let $A$ be an adjacency matrix of $(G,\alpha)$, which means that for any two vertices $u,v$ of $G$ the $uv$-th entry of $A$ is $\alpha(e)$ if $e=uv$ is an edge of $G$, and $0$ otherwise.
\smallskip

If $A$ is non-singular, the {\em inverse} of $(G,\alpha)$ is the edge-labeled graph $(H,\beta)$ determined by the adjacency matrix equal to the inverse $A^{-1}$ of $A$. We thus assume that $G$ and $H$ share the same vertex set, and $e=uv$ is an edge of $H$ if the $uv$-th entry of $A^{-1}$ is non-zero, in which case this $uv$-th entry is the label $\beta(e)$ of $e$. Obviously, the inverse defined this way is unique up to graph isomorphism preserving edge-labels.
\smallskip

Inverses of edge-labeled graphs as introduced above were studied, for example, in \cite{Bapat,PP,Pavl,Ye}. Edge-labels in \cite{Bapat,Pavl} were even allowed to be elements of a not necessarily commutative ring, and a formula for an inverse graph to $(G,\alpha)$ was given in both papers in the special case of a bipartite graph $G$ with a unique perfect matching and with multiplicatively invertible $\alpha$-labels on matched edges.
\smallskip

Returning to real-valued labelings, one may further restrict to graphs $(G,\alpha)$ in which all edge-labels are positive. In such a case of a {\em positively labeled} graph with a non-singular adjacency matrix, its inverse will typically contain both positive as well as negative entries. We will then say that a positively labeled graph $(G,\alpha)$ with a non-singular adjacency matrix $A$ is {\em positively invertible} if $A^{-1}$ is diagonally similar (signable, in the terminology of \cite{AK}) to a non-negative matrix. Positively invertible graphs with integral edge-labels have been studied in detail in \cite{BNP,AK,McMc} and the last paper also contains a nice survey of development in the study of inverses of graphs.
\smallskip

All this work, however, was initiated by the influential paper \cite{Gods} about inverses of trees with all edge-labels equal to $1$, extending earlier observations of \cite{CGS}. By \cite{Har}, an adjacency matrix of a tree is invertible if and only if the tree has a (unique) perfect matching. In terms of our definition, each such tree is automatically invertible. The much stronger result of \cite{Gods} says that every tree with a perfect matching is {\em positively invertible} and its positive inverse is a {\em simple graph} (containing no loops) with every edge carrying the {\em unit label} again. A formula for determining the inverse of a tree with a perfect matching in terms of alternating paths appeared later in \cite{PK}, which was afterwards extended to bipartite graphs with a unique perfect matching in \cite{BNP,AK,Bapat,McMc,Pavl}. For completeness, graphs arising as inverses of trees with a perfect matching, and as inverses of bipartite graphs with a unique perfect matching that remain bipartite after contracting the matching, were characterized in \cite{NP} and \cite{PP}, respectively, and self-inverse graphs in the latter family were classified in \cite{SC}.
\smallskip

In this situation it is natural to ask what one can do in the case of edge-labeled graphs with a {\em singular} adjacency matrix. An equally natural move is to consider `inverting' the matrix by taking one of the generalizations of matrix inverses, such as the Moore-Penrose inverse, or the Drazin inverse, or a special case of the latter known as the group inverse. In the instance of a (square) symmetric matrix $A$ all these inverses coincide and we will call the corresponding matrix simply a {\em generalized inverse} and denote it by $A^*$. Its definition in the symmetric case is simple. Since a real symmetric $n\times n$ matrix $A$ is orthogonally diagonalizable, there is an orthogonal matrix $P$ such that $PAP^T=D$, where $D={\rm diag}(\lambda_1,\ldots,\lambda_k,0,\ldots,0)$ is the diagonal matrix of eigenvalues of $A$, with $k={\rm rank}(A)$ non-zero eigenvalues $\lambda_1,\ldots,\lambda_k$. Letting $D^*={\rm diag} (\lambda_1^{-1}, \ldots,\lambda_k^{-1},0,\ldots,0)$, the {\em generalized inverse} $A^*$ of $A$ is simply given by $A^*=PD^*P^T$, that is, both $A$ and $A^*$ are conjugate to their corresponding diagonal matrices by the {\em same} orthogonal matrix $P$. Note that $A^*$ is again symmetric, and $A^*$ coincides with $A^{-1}$ if $A$ is non-singular.
\smallskip

Motivated by this, we define the {\em generalized inverse} of an edge-labeled graph $(G,\alpha)$ with adjacency matrix $A$ to be the labeled graph $(G^*,\alpha^*)$ with adjacency matrix $A^*$, the generalized inverse of $A$. As before, $G$ and $G^*$ are assumed to have the same vertex set, and $e=uv$ is an edge of $G^*$ if and only if the $uv$-th entry of $A^*$ is non-zero, and then this entry is also the label $\alpha^*(e)$ of $e$. And, again, note that $G^*$ is well defined up to isomorphism preserving edge labels.
\smallskip

Observe that this way of defining generalized inverses of edge-labeled graphs is in line with the original motivation of considering graph inverses which comes from chemistry. Namely, there appear to be fewer methods for estimating the smallest positive eigenvalue of a graph in contrast to a larger number of techniques for bounding the largest positive eigenvalue. For graphs representing structure of molecules, however, the smallest positive eigenvalue is a meaningful parameter in quantum chemistry. If such a graph has an inverse, one may hope to increase the number of techniques for estimating its smallest positive eigenvalue by passing to bounds on the largest positive eigenvalue of the inverse graph. This feature remains present also for our generalized inverses.
\smallskip

The main result of this paper is a formula for calculating the generalized inverse of an arbitrary tree $T$ equipped with the constant label $\alpha(e)=1$ on every edge $e$; we will simply refer to $T$ and omit $\alpha$ in the notation. Also, we will assume throughout that $T$ is non-trivial, that is, containing at least two vertices.
\smallskip

To state the result we need to introduce a few concepts. If $M$ is a matching in $T$ and $u,v$ are distinct vertices of $T$, an $u_Mv$ {\em alternating path} is a $u-v$ path $P$ in $T$ whose edges belong alternately to $M$ and not to $M$, with the condition that the first and last edges of $P$ (that is, those incident to $u$ and $v$) both belong to $M$. Note that for given $u,v$ and $M$ such an $u_Mv$ alternating path in $T$ may not exist (it never does if the distance between $u$ and $v$ is even) but if it does, it is unique (meaning that $u_Mv$ is the same for possibly different maximum matchings $M$). Note that, in particular, any edge $xy$ of $M$ is an $x_My$ alternating path itself. For a pair of distinct vertices $u,v$ in $T$ at an odd distance we let $\sigma(u,v)$ be equal to $+1$ or $-1$ depending on whether their distance is congruent to $+1$ or $-1$ mod $4$. Further, for such a pair we let $\mu_T(u,v)$ be the $\sigma(u,v)$-multiple of the number of all the {\em maximum} matchings $M$ in $T$ for which there exists an $u_Mv$ alternating path in $T$, and we let $\mu_T(u,v)=0$ otherwise. Finally, we let $m(T)$ denote the total number of maximum matchings of $T$. In this terminology and notation our main result is:

\bt\label{main}
Let $T$ be a tree with vertex set $V$. Then, its generalized inverse $(T^*,\alpha^*)$ has two distinct vertices $u,v\in V$ joined by an edge $e$ if and only if the original tree $T$ contains an $u_Mv$ alternating path for some maximum matching $M$ of $T$, and the label of $e$ is then given by $\alpha^*(e)=\alpha^*(uv) = \mu_T(u,v)/m(T)$.
\et

This result generalizes the original findings of \cite{PK} on inverses of trees with a (unique) perfect matching, which are those with $m(T)=1$.

\section{The strategy}

The generalized inverse $A^*$ of a symmetric matrix $A$ was introduced by way of conjugating the associated diagonal matrices by the {\em same} matrix, denoted $P$ in the Introduction. By elementary linear algebra, however, all such matrices $P$ represent bases of the corresponding eigenspaces; we summarize this in form of an observation that needs no proof.

\bl\label{lem:01}
Two symmetric square matrices $A$ and $B$ of the same dimension are generalized inverses of each other if and only if they have the same null-space and, for every non-zero eigenvalue $\lambda$ of $A$ the quantity $\lambda^{-1}$ is an eigenvalue of $B$ and the corresponding eigenspaces of $A$ and $B$ are identical. \hfill $\Box$
\el
\smallskip

It turns out that Lemma \ref{lem:01} can be given an equivalent form with no explicit reference to actual values of $\lambda\ne 0$. For a positive integer $n$ we let $[n]=\{1,2,\ldots,n\}$.

\bp\label{prop:gener}
Let $A=(a_{ij})$ and $B=(b_{ij})$ be symmetric $n\times n$ matrices. Then, $B$ is the generalized inverse of $A$  if and only if $A$ and $B$ have the same null-spaces, and every eigenvector $f:\ [n]\to {\mathbb R}$ of $A$ corresponding to a non-zero eigenvalue of $A$ satisfies
\be\label{eq:01}
\sum_{j\in [n]} b_{ij} \sum_{k\in [n]} a_{jk} f(k) = f(i)\ \ {\rm for\ every\ } i\in [n]\ .
\ee
\ep

\pr \
Let $f:\ [n]\to {\mathbb R}$ be an eigenvector of $A$ corresponding to a non-zero eigenvalue $\lambda$ of $A$; that is,
\be\label{eq:02}
\sum_{k\in [n]} a_{jk}f(k)=\lambda f(j) \ \ {\rm for\ every\ } j\in [n]\ .
\ee
Multiplying both sides of (\ref{eq:02}) by $b_{ij}$ and summing over all $j\in [n]$ gives
\be\label{eq:03}
\sum_{j\in [n]} b_{ij} \sum_{k\in [n]} a_{jk} f(k) = \lambda \left( \sum_{j\in [n]} b_{ij}f(j) \right) \ .
\ee
If $B=A^*$, then $f$ is also an eigenvector of $B$ corresponding to $\lambda^{-1}$ by Lemma \ref{lem:01}. This means that $\sum_{j\in [n]} b_{ij}f(j)=\lambda^{-1}f(i)$, and hence (\ref{eq:03}) reduces to (\ref{eq:01}). Of course, if $B=A^*$, then both matrices have the same null-space. This proves necessity.
\smallskip

To prove sufficiency, assuming (\ref{eq:01}) the left-hand side of (\ref{eq:03}) reduces to $f(i)$ and a subsequent  division by $\lambda\ne 0$ gives $\lambda^{-1}f(i) = \sum_{j\in [n]} b_{ij}f(j)$. This shows that $f$ is, at the same time, an eigenvector of $B$ for the eigenvalue $\lambda^{-1}$. Combining this with the assumption on null-spaces we conclude that $B=A^*$ by Lemma \ref{lem:01}.
\ebox

Proposition \ref{prop:gener} may also be of general interest by giving quite straightforward necessary and sufficient conditions for generalized invertibility, suitable for particular applications.
\smallskip

With this tool in hand we will now outline the strategy of proving our main result. In the terminology and notation introduced before the statement of Theorem \ref{main}, for every tree $T$ with an adjacency matrix $A$ we need to show that if $B$ is a matrix (indexed the same way as $A$) with  $uv$-th entry equal to $\mu_T(u,v)/m(T)$, then $B=A^*$. In the light of Proposition \ref{prop:gener} our task reduces to showing validity of the following statement, where for any vertex $x$ of $T$ we let $T(x)$ denote the set of neighbours of $x$ in $T$.

\bp\label{prop:eigen}
Let $A$ be an adjacency matrix of a tree $T$ on a vertex set $V$. Assume that the following two conditions {\rm (a), (b)} are fulfilled:
\bi
\item[{\rm (a)}] Every eigenvector $f:\ V\to {\mathbb R}$ of $A$ corresponding to a non-zero eigenvalue satisfies
\be\label{eq:1}
\sum_{v\in V} \mu_T(u,v) \sum_{w\in T(v)} f(w) = m(T)\cdot f(u)\ \ {\rm for\ every\ } u\in V\ .
\ee
\item[{\rm (b)}] For very eigenvector $g:\ V\to {\mathbb R}$ of $A$ corresponding to a zero eigenvalue of $A$,
\be\label{eq:2}
\sum_{v\in T(u)}g(v)=0 \ \ {\rm implies} \ \ \sum_{v\in V}\mu_T(u,v)g(v)=0\ \ {\rm for\ every\ } u\in V\ .
\ee
\ei
Then, the square matrix $B$ (indexed as $A$) with $uv$-th entry equal to $\mu_T(u,v)/m(T)$ for every $u,v\in V$ is the generalized inverse of $A$.
\ep

\pr \
The statement is essentially the sufficient condition of Proposition \ref{prop:gener}. The set $[n]$ of indices has been replaced here by $V$ and the values equal to $1$ in the $(0,1)$-adjacency matrix $A$ of $T$ correspond to entries appearing in summations over neighborhoods of vertices of $T$. To obtain the condition (a) one just has to multiply both sides of (\ref{eq:01}) by $m(T)$ in the new notation; the condition (b) is the null-spaces assumption.
\ebox

This shows that to prove Theorem \ref{main} it is sufficient to prove validity of conditions (a) and (b) from Proposition \ref{prop:eigen}. We divide this task into three parts. In Section \ref{sec:aux} we prove auxiliary results regarding maximum matchings in trees, which may be of independent interest due to their generality. Parts (a) and (b) of Proposition \ref{prop:eigen} are proved in Sections \ref{sec:nonz} and \ref{sec:zero}, again by results that may be of broader interest. The final Section \ref{sec:fin} wraps up the proof of Theorem \ref{main} and contains a few concluding remarks.

\section{Auxiliary results on maximum matchings}\label{sec:aux}

Aiming for establishing the condition (a) of Proposition \ref{prop:eigen} we first prove two auxiliary results regarding numbers of maximum matchings satisfying some extra conditions. We will again refer to the terminology and notation introduced earlier. In addition, we will say that an $u{-}v$ path $P$ in a tree is {\em matched} if it has the form $u_Mv$ for some matching $M$ of the tree. In such a situation we will also say that $M$ turns $P$ into a matched path, or that $M$ makes $P$ a matched path. Since we will not need to evaluate the function $\mu$ for different trees we will omit the subscript and write just $\mu(u,v)$ instead of $\mu_T(u,v)$. To state our first result, we will say that a vertex $u$ of a tree $T$ is {\em saturated} if $u$ belongs to {\em every} maximum matching of $T$.

\bl\label{lem:1}
Let $u$ be a saturated vertex of $T$. Then, for every vertex $w$ of $T$ we have
\be\label{eq:41}
\sum_{v\in T(w)} \mu(u,v) \ = \
\begin{cases}
\ \ m(T) & {\rm if\ }  w = u\ , \\
\ \ 0    & {\rm if\ } w \ne u\ .
\end{cases}
\ee
\el

\pr \
Note that the assertion is vacuously true if $w$ is at an odd distance from $u$. Assume thus that the distance between $w$ and $u$ is even. If $w=u$, then $\mu(u,v)$ for $v\in T(w)=T(u)$ is simply the number of maximum matchings in $T$ containing the edge $uv$. The sum in (\ref{eq:41}) over all neighbours $v$ of $u$ is then equal to the number of maximum matchings in $T$ that cover the vertex $u$, which is equal to $m(T)$ since $u$ is assumed to be saturated.
\smallskip

If $w\ne u$ has even distance from $u$, let $u\ldots yw$ be the unique $u{-}w$ path in $T$ and let $Q$ be its $u{-}y$ subpath. Let $m(Q)$ be the number of maximum matchings in $T$ that turn $Q$ into a matched path. Further, for every neighbour $v$ of $w$ distinct from $y$ let $m(Q,wv)$ be the number of maximum matchings as above but containing the edge $wv$ as well. The contribution of the neighbour $v=y$ of $u$ to the sum in (\ref{eq:41}) is $\mu(u,y) = \varepsilon m(Q)$ for some $\varepsilon = \pm 1$ depending on the (odd) length of $Q$. The contribution of all the remaining neighbours $v\ne y$ of $w$ is, due to different residue class of the paths $u{-}ywv$ compared with $Q$, equal to the sum $-\varepsilon \sum m(Q,wv)$ ranging over all $v\in T(w){\setminus}\{y\}$. Thus, (\ref{eq:41}) will be established if we show that, for every $w\ne u$ at an even distance from $u$,
\be\label{eq:42}
\sum_{v\in T(w){\setminus}\{y\}} m(Q,wv) = m(Q)\ .
\ee
Suppose that (\ref{eq:42}) is invalid. It is obvious that this happens if and only if there is a maximum matching $M$ in $T$ that turns $Q$ into a matched path but does not cover $w$. But then, trading the matched edges on the path $Q\cup\{yw\}$ for the unmatched ones we would obtain from $M$ a maximum matching in $T$ that does not cover the vertex $u$. This, however contradicts the saturation assumption on $u$, and completes our proof.
\ebox

To be able to formulate our second result in this section, for every vertex $x$ of our tree $T$ we let $T_2(x)$ denote the multi-set consisting of end-vertices $v$ of all walks of length $2$ of the form $xwv$ that start at $x$. Note that if $x$ has valency $d$, then the vertex $x$ itself appears $d$ times in $T_2(x)$, due to the $d$ walks of the form $xwx$ for $w\in T(x)$.
\smallskip

\bl\label{lem:2}
Let $u$ and $x$ be arbitrary vertices of a tree $T$. Then,
\be\label{eq:21}
\sum_{v\in T_2(x)}\mu(u,v) \ = \
\begin{cases}
\ \ m(T) & {\rm if\ }  x\in T(u)\ , \\
\ \ 0    & {\rm if\ } x\notin T(u)\ .
\end{cases}
\ee
\el

\pr \
Observe that the statement is vacuously true if $x$ and $u$ are at an even distance in $T$, in which case the result of (\ref{eq:21}) is zero because all the values $\mu(u,v)$ are automatically equal to $0$. We therefore assume that $x$ is at an odd distance from $u$ to allow for matchable $u{-}v$ paths in the sum appearing in (\ref{eq:21}).
\smallskip

Let $P$ be the unique $u{-}x$ path in $T$; by our assumptions $P$ has odd length. We may assume that either $P$ has the form $u\ldots yzx$, with $yzx$ being a sub-path of $P$ of length $2$ and $y\ne u$, or $P$ is the single edge $zx$ with $z=u$; this choice of notation will be handy later. Let $T_2(x;\neg\,  z)$ and $T_2(x;z)$ be the sets of end-vertices $v\ne x$ of {\em paths} $xwv$ for $w\ne z$ and $w=z$, respectively; here and henceforth $z=u$ if $P=ux$. Letting $d(x)$ be the valency of $x$ in $T$, the multi-set $T_2(x)$ consists of $d(x)$ entries equal to $x$ and of single entries in $T_2(x;\neg\, z)\cup T_2(x;z)$. We proceed with a separate evaluation of contributions of the $d(x)$ entries $v=x$, the entries $v$ from the set $T_2(x;\neg\,  z)$, and the entries $v$ from $T_2(x;z)$, to the sum in (\ref{eq:21}).
\smallskip

{\em Contribution of the $d(x)$ entries $v=x$}. \ Let $m(P)$ be the set of all maximum matchings in $T$ that turn $P$ into a matched path. Then, for each of the $d(x)$ neighbours $w\in T(x)$ the walk $xwv$ for $v=x$ contributes to the sum in (\ref{eq:21}) by $\mu(u,x) = \varepsilon m(P)$ for some $\varepsilon = \pm 1$ depending on the congruence class of the (odd) length of $P$ mod $4$. The total contribution to (\ref{eq:21}) of such entries is therefore equal to $\varepsilon d(x) m(P)$.
\smallskip

{\em Contribution of $v\in T_2(x;\neg\, z)$}. \ Let $W=T(x){\setminus}\{z\}$, the set of all neighbours of $x$ distinct from $z$. For every $w\in W$  we let $m(P,w)$ denote the number of maximum matchings in $T$ turning $P$ into a matched path and covering the vertex $w$. Each such matching induces an alternating path $u\ldots xwv$ for some $w\in W$ and $v\in T_2(x;\neg\, z)$. Note that the (odd) length of each such alternating path is in a different congruence class mod $4$ compared with $P$. The total contribution of all the $\mu(u,v)$ terms to (\ref{eq:21}) taken over all vertices $v\in T_2(x;\neg\, z)$ is thus equal to
\[ \sum_{w\in W}\sum_{v\in T(w){\setminus}\{z\}} \mu(u,v) = -\varepsilon \sum_{w\in W}m(P,w)\ .\]
For contribution of $v\in T_2(x;z)$, we need to consider the cases $x\in T(u)$ and $x\notin T(u)$ separately.
\smallskip

{\em Contribution of $v\in T_2(x;z)$ if $x\in T(u)$.} \ Here the set $T_2(x;z)$ for $z=u$ coincides with the set of neighbours $v$ of $u$ distinct from $x$. For every such vertex $v\in T(u){\setminus}\{x\}$ its contribution to the sum in (\ref{eq:21}) is $\mu(u,v)$, which is the number of maximum matchings in $T$ containing the edge $uv$. Letting $m(u)$ be the number of maximum matchings in $T$ that cover $u$, it follows that if $P=ux$, the net contribution to (\ref{eq:21}) in this case is equal to $m(u) - m(P)$.
\smallskip

{\em Contribution of $v\in T_2(x;z)$ if $x\notin T(u)$.} \ Let $Q$ be the $u{-}y$ sub-path of $P$ and let $m(Q)$ and $m(Q,z,\neg\,  zx)$ be the number of maximum matchings in $T$ that turn $Q$ into a matched path, and the number of such paths that in addition cover $z$ but do not contain the edge $zx$, respectively. The contribution of all vertices $v\in T_2(x;z) {\setminus}\{y\}$ via matched $u{-}yzv$ paths (of the same length as $P$) to the sum in (\ref{eq:21}) is clearly equal to $\varepsilon m(Q,z,\neg\, zx)$. It remains to consider $v=y$; since the length of the $u{-}y$ path $Q$ has different odd parity mod $4$ compared with $P$, the vertex $v=y$ contributes to (\ref{eq:21}) by $\mu(u,y) = -\varepsilon m(Q)$. It follows that for $x\notin T(u)$ the net contribution of this kind is $\varepsilon [m(Q,z,\neg\,  zx) - m(Q)]$. To express this in a more suitable form, observe that if we let $m(Q,\neg\,  z)$ denote the number of maximum matchings in $T$ turning $Q$ into a matched path and avoiding $z$, then we have $m(Q)=m(P) + m(Q,z, \neg\, zx) + m(Q,\neg\, z)$. This enables us to express the contribution to (\ref{eq:21}) in this last case as $-\varepsilon[m(Q,\neg\, z) + m(P)]$.
\smallskip

We proceed with writing the left-hand side of (\ref{eq:21}) as the sum of the three types of contributions that we have just determined. Using $d(x)=|W|+1$ and replacing $|W|$-multiples by sums $\sum_{w\in W}$ of equal terms the left-hand side of (\ref{eq:21}) evaluates as follows:
\be\label{eq:22}
\sum_{v\in T_2(x)}\mu(u,v) \ = \
\begin{cases}
\ \ \ \ \ \sum_{w\in W} \left( m(P) - m(P,w) \right) + m(u)\  & {\rm if\ }  x\in T(u)\ , \\
\ \ \varepsilon\left[ \sum_{w\in W} \left( m(P)-m(P,w) \right) - m(Q,\neg\,  z) \right]  & {\rm if\ } x\notin T(u)\ .
\end{cases}
\ee

Comparing (\ref{eq:21}) with (\ref{eq:22}) it is clear that to prove Lemma \ref{lem:2} it is sufficient to show that \be\label{eq:23}
\sum_{w\in W} \left( m(P) - m(P,w)  \right) \ = \
\begin{cases}
\ \ m(T) - m(u) & {\rm if\ } x\in T(u)\ , \\
\ \  m(Q,\neg\,  z)\ & {\rm if\ } x\notin T(u)\ .
\end{cases}
\ee

To establish (\ref{eq:23}), we will consider again the $u{-}y$ path $Q$ for $x\notin T(u)$ but we also let $Q=\emptyset$ if $x\in T(u)$ (and then we again use $z=u$). For each $w\in T(x){\setminus}\{z\}$ let ${\cal M}(Q,xw,\neg\,  z)$ and ${\cal M}(Q,zx,\neg\,  w)$ denote the set of maximum matchings of $T$ turning $Q$ into a matched path (this condition is vacuous if $Q=\emptyset$) and, respectively, containing the edge $xw$ while not covering $z$, and those containing the edge $zx$ while not covering $w$. For each $M\in {\cal M}(Q,xw,\neg\,  z)$ the matching $M'= M{\setminus}\{xw\} \cup\{zx\}$ obviously belongs to ${\cal M}(Q,zx,\neg\,  w)$, and the mapping $M\mapsto M'$ is a bijection ${\cal M}(Q,xw,\neg\,  z) \to  {\cal M}(Q,zx,\neg\,  w)$. It follows that
\be\label{eq:24}
\sum_{w\in W} |{\cal M}(Q,zx,\neg\,  w)| = \sum_{w\in W} |{\cal M}(Q,xw,\neg\,  z)| \ .
\ee

The left-hand side of (\ref{eq:24}) is obviously equal to the left-hand side of (\ref{eq:23}), irrespective of the position of $x$ relative to $T(u)$. If $x\in T(u)$, then $z=u$ and the right-hand side of (\ref{eq:24}) is clearly equal to $m(T) - m(u)$. Finally, if $x\notin T(u)$, the sum on the right-hand side of (\ref{eq:24}) is equal to $m(Q,\neg\,  z)$ since every maximum matching that turns $Q$ into a matched path but avoids $z$ must contain an edge $xw$ for some $w\in W$. This proves (\ref{eq:23}), and hence completes the proof of Lemma \ref{lem:2}.
\ebox

\section{Eigenvectors for non-zero eigenvalues}\label{sec:nonz}

In this section we establish validity of the statement in part (a) of Proposition \ref{prop:eigen} and we begin with saturated vertices $u$ of $T$. Actually, for saturated vertices we prove the following statement which is much more general than the one of part (a) in Proposition \ref{prop:eigen} in that it is valid for {\em arbitrary} real functions defined on the vertex set of a tree, and not just for eigenvectors corresponding to a non-zero eigenvalue.

\bp\label{prop:satur}
Let $T$ be a tree with vertex set $V$ and let $u\in V$ be a saturated vertex. Then, every function $f: \ V\to {\mathbb R}$ satisfies
\be\label{eq:4}
\sum_{v\in V} \mu_T(u,v) \sum_{w\in T(v)} f(w) = m(T)f(u) \ .
\ee
\ep

\pr \
Changing the order of summation in (\ref{eq:4}) turns this equation into the following equivalent form:
\be\label{eq:4.2}
\sum_{w\in V} c_u(w)f(w) = m(T)f(u)\ , \ \  {\rm where} \ \ \ c_u(w) = \sum_{v\in T(w)} \mu(u,v) \ .
\ee
But in Lemma \ref{lem:1} we have shown that $c_u(w) = m(T)$ if $w=u$, and $c_u(w)=0$ for every $w\in V{\setminus}\{u\}$, which makes (\ref{eq:4.2}) automatic.
\ebox

We state and prove yet another auxiliary result that may also be of independent interest since it is again valid for all real functions defined on vertices of a tree; moreover, we will use it to handle the part (a) of Proposition \ref{prop:eigen} for unsaturated vertices.

\bp\label{prop:complex}
Let $T$ be a tree with vertex set $V$ and let $u\in V$ be an arbitrary vertex. Then, for every function $f: \ V\to {\mathbb R}$ we have
\be\label{eq:5.1}
\sum_{v\in V} \mu_T(u,v) \sum_{w\in T(v)} \sum_{x\in T(w)} f(x) = m(T)\sum_{y\in T(u)}f(y) \ .
\ee
\ep

\pr \
Recall the notation introduced in section \ref{sec:aux}, where we denoted by $T_2(x)$ the multi-set of end-vertices $v$ of all walks of length $2$ starting at $x$, that is, walks of the form $xwv$. We also reiterate that if $x$ has valency $d$, then $x$ itself appears $d$ times in $T_2(x)$, which is because of the existence of $d$ walks of the form $xwx$ for $w\in T(x)$. With this notation it is easy to see that, by changing the order of summation, one can manipulate the equation (\ref{eq:5.1}) into the equivalent form
\be\label{eq:5.2}
\sum_{x\in V} C_u(x) f(x) = m(T)\sum_{y\in T(u)}f(y)\ ,\ \ {\rm where}\ \ \ C_u(x)= \sum_{v\in T_2(x)}\mu_T(u,v)\ .
\ee
Again, most of the work needed has already been done. Namely, in Lemma \ref{lem:2} we have proved that $C_u(x) = m(T)$ if $x\in T(u)$, while $C_u(x)=0$ if $x\notin T(u)$. Validity of (\ref{eq:5.2}) is now obvious.
\ebox

We are now ready to tackle the part (a) of Proposition \ref{prop:eigen} in its full generality. It is interesting to note that this is the only place where we are going to use the assumption that the real functions we work with are actually eigenvectors corresponding to non-zero eigenvalues of a tree. Because of Proposition \ref{prop:satur} we need to do this only for non-saturated vertices but, as we shall see, most of the arguments will be valid in general.

\bp\label{prop:unsat}
Let $T$ be a tree with vertex set $V$ and let $f:\ V\to {\mathbb R}$ be an eigenvector of an adjacency matrix of $T$ corresponding to a non-zero eigenvalue. Then, for every vertex $u$ of $T$,
\be\label{eq:3'}
\sum_{v\in V} \mu_T(u,v) \sum_{w\in T(v)} f(w) = m(T)f(u) \ .
\ee
\ep

\pr
As alluded to, we have already established (\ref{eq:3'}) for saturated vertices $u$ in Proposition \ref{prop:satur}, but in what follows we will not make use of this assumption. Let $f$ be an eigenvector as in the statement, corresponding to a non-zero eigenvalue $\lambda$. Take the equation (\ref{eq:5.1}) from Proposition \ref{prop:complex} and substitute there $\lambda f(u)$ for $\sum_{y\in T(u)}f(y)$ and $\lambda f(w)$ for $\sum_{x\in T(w)} f(x)$, the latter for every $w\in T(v)$. After cancelling the non-zero term $\lambda$ on both sides we clearly obtain (\ref{eq:3'}), completing the proof.
\ebox

It should now be clear that a proof of part (a) of Proposition \ref{prop:eigen} can be obtained solely on the basis of Lemma \ref{lem:2} and Proposition \ref{prop:unsat} (by incorporating Proposition \ref{prop:complex}). We have expanded the presentation slightly by including Proposition \ref{prop:satur} (based on Lemma \ref{lem:1}) and Proposition \ref{prop:complex} because of possible independent interest due to their validity for {\em all} real functions defined on the vertex set of a tree.

\section{Eigenvectors for the zero eigenvalue}\label{sec:zero}

Part (b) of Proposition \ref{prop:eigen} simply says that an adjacency matrix $A$ of a tree $T$ and the matrix $B$ (indexed the same as $A$) with entries $\mu_T(u,v)/m(T)$ have the same null-spaces. By elementary linear algebra this is equivalent to $B$ being convertible to $A$ by a sequence of elementary row operations, or, $B$ being {\em Gaussian equivalent} to $A$, for short. Since the matrix $C=m(T){\cdot}B$ is Gaussian equivalent to $B$ and contains integer entries $\mu_T(u,v)$ only, it will be of advantage to focus on Gaussian equivalence of $C$ with $A$, which we will do next.

\bp\label{prop:Gauss}
Let $A$ be an adjacency matrix of a tree $T$ with vertex set $V$ and let $C$ be a matrix (indexed the same way as $A$ by elements of $V$) with entries $C_{uv}=\mu_T(u,v)$ for every $u,v\in V$. Then, $A$ and $C$ are Gaussian equivalent over the field of rational numbers.
\ep

\pr
We will use induction on $|V|$, observing that the statement is trivial for stars. Let $T$ be a non-star tree on a vertex set $V$, with symmetric adjacency matrix $A$ indexed by elements of $V$. Let $(V_1,V_2)$ be the bipartition of $V$ with parts $V_1$ and $V_2$. Since no edge of $T$ is incident to a pair of vertices from the same part, the $uv$-entries of $A$ are zero whenever $u,v\in V_1$ or $u,v\in V_2$. We will assume that $C$ is indexed the same way as $A$, and as $\mu_T(u,v)$ can only be non-zero if $u$ and $v$ belong to different parts, we again have $uv$-th entries of $C$ equal to zero if $u,v\in V_1$ or $u,v\in V_2$. In what follows we will set the stage for the induction step.
\smallskip

Let $x$ be a pendant vertex of a $T$ such that $x$ is the end-vertex of a longest path in $T$; we may without loss of generality assume that $x\in V_1$. Let $y\in V_2$ be the vertex adjacent to $x$ in $T$; clearly, $y$ must be saturated in $T$. We now change a bit the notation for sub-trees of $T$ used in the previous section to avoid confusion. By $T{\setminus}x$ we denote the subtree of $T$ obtained by deleting the vertex $x$, and by and $T{\setminus}xy$ the subtree of $T{\setminus}x$ arising as the only non-trivial connected component that remains after deleting $y$ from $T{\setminus}x$. Let $z\in V_1$ be the (uniquely determined) non-pendant vertex of $T$ adjacent to $y$; note that $z$ is a vertex of $T{\setminus}xy$. To facilitate further explanation, we let $|V|=n$ and we will (without loss of generality) assume that $x$, $y$ and $z$ are the indices of the $n$-th, $(n-1)$-st and $(n-2)$-nd row and column of the (symmetric) matrices $A$ and $C$.
\smallskip

Let $A_x$ and $A_{xy}$ be adjacency matrices of $T{\setminus}x$ and $T{\setminus}xy$, obtained, respectively, by deleting the last row and column from $A$ and the last row and column from $A_x$. We also introduce the corresponding matrix $C_x$ with $uv$-th entry for $u\in V_1{\setminus}\{x\}$ and $v\in V_2$ equal to $\mu_{T{\setminus}x}(u,v)$, and $C_{xy}$ with $uv$-th entry for $u\in V_1{\setminus}\{x\}$ and $v\in V_2{\setminus}\{y\}$ equal to $\mu_{T{\setminus}xy}(u,v)$; indexing of the two matrices is assumed to be inherited from $A_x$ and $A_{xy}$.
\smallskip

Take now an arbitrary vertex $u\in V$ and fix it. We begin by assuming that $u\in V_1$ (note that, having placed our pendant vertex $x$ in $V_1$, the case when $u\in V_2$ needs a separate consideration). Observing that a maximum matching $M$ in $T$ either contains the vertex $x$ and hence the edge $xy$ (and such matchings are in a one-to-one correspondence with maximum matchings in $T{\setminus}xy$), or $M$ does not contain $x$ and then $M$ must still contain $y$; such matchings are exactly the maximum matchings in $T{\setminus}x$ that contain $y$.
\smallskip

In fact, if the valency $d(y)$ of $y$ in $T$ is at least $3$, then $y$ is saturated also in $T{\setminus}x$. Indeed, by the choice of $x$ as an end-vertex of a longest path in $T$, if $d(y)\ge 3$ then $y$ is adjacent to some {\em pendant} vertex $x'\notin \{x,y\}$ in $T{\setminus}x$ and hence $y$ is saturated in $T{\setminus}x$. But if $d(y) = 2$, the vertex $y$ may or may not be saturated in $T{\setminus}x$, and if it is not, then no maximum matching $M'$ of $T{\setminus}x$ containing $y$ is extendable to a maximum matching of $T$. But then so such matching $M'$ can  contribute to the value of $\mu_T(u,v)$ for any vertices $u,v$ of $T{\setminus}x$. Because of this subtlety we call the vertex $y$ {\em exceptional} if $d(y) = 2$ and $y$ is not saturated in $T{\setminus}x$. Letting $\theta=0$ if $y$ is exceptional and $\theta=1$ otherwise, this analysis implies that
\[\begin{array}{c}
\mu_T(u,v) = \theta{\cdot}\mu_{T{\setminus}x}(u,v) + \mu_{T{\setminus}xy}(u,v) \ \ {\rm for\ every}\ v\in V_2\ {\rm if}\ u\ne x, \ \ {\rm and} \\
\mu_T(x,y)=m(T{\setminus}{xy})\ \ {\rm and} \ \ \mu_T(x,v) = -\mu_{T{\setminus}xy}(z,v)\ \ {\rm for} \ \ v\in V_2{\setminus}\{y\}\ .
\end{array}\]
It follows that rows $C(u)$ of $C$ indexed by elements $u\in V_1{\setminus}\{x\}$ and by $x$ can be conveniently displayed as follows, with zeros in columns corresponding to vertices in $V_1{\setminus}\{x\}$:
\begin{equation}\label{eq:10}
\bordermatrix{
C & & {\rm column}\ v\in V_2{\setminus}\{y\} & & {\rm column}\ y \ & \ {\rm column}\ x & \cr
{\rm row}\ C(u)\ & \ \ldots & \theta\mu_{T{\setminus}x}(u,v)+\mu_{T{\setminus}xy}(u,v) & \ldots & \mu_{T{\setminus}x}(u,y) & 0  \cr
{\rm row}\ C(x)\ & \ \ldots & -\mu_{T{\setminus}xy}(z,v) & \ldots & \mu_T(x,y) & 0 \cr }
\end{equation}
We will show that, for every $u\in V_1{\setminus}\{x\}$, the row $C(u)$ in (\ref{eq:10}) is a linear combination of rows of $A$. Indeed, for such a $u$ consider the row $C_x(u)$ of the matrix $C_x$ which has the form
\begin{equation}\label{eq:11}
\bordermatrix{
C_x & &  {\rm column}\ v\in V_2{\setminus}\{y\} & & {\rm column}\ y \cr
{\rm row}\ C_x(u)\  & \ \ldots & \mu_{T{\setminus}x}(u,v) & \ldots & \mu_{T{\setminus}x}(u,y) \ \cr }
\end{equation}
By the induction hypothesis applied to the tree $T{\setminus}x$, the row (\ref{eq:11}) is a linear combination of rows $A_x(w)$ of $A_x$, that is, $C_x(u)=\sum_{w\in V_1{\setminus}\{x\}}\delta_wA_x(w)$ for some rational numbers $\delta_w$; the restriction in the sum to elements in $V_1{\setminus}\{x\}$ comes from the fact that rows of $A_x$ indexed by elements of $V_2$ have zeros in columns corresponding to $V_2$. Note that rows $A_x(w)$ for $w\in V_1{\setminus}\{x\}$, of length $n-1$, differ from those of $A(w)$, of length $n$, only in the $n$-th coordinate (indexed by $x$) of $A(w)$, equal to zero. It follows that, taking the rows $A(w)$ of the matrix $A$ instead of the rows $A_x(w)$ and using the {\em same} coefficients $\delta_w$ as above, the resulting linear combination $S_x(u)=\sum_{w\in V_1{\setminus}\{x\}} \delta_wA(w)$ gives
\begin{equation}\label{eq:12}
\bordermatrix{
  & & {\rm column}\ v\in V_2{\setminus}\{y\} & & {\rm column}\ y \ & \ {\rm column}\ x \ \cr
S_x(u)\ & \ \ldots & \mu_{T{\setminus}x}(u,v) & \ldots & \mu_{T{\setminus}x}(u,y) & 0 \ \cr }
\end{equation}
Next, for $u\in V_1{\setminus}\{x\}$ consider the $u$-th row $C_{xy}(u)$ of the matrix $C_{xy}$:
\begin{equation}\label{eq:13}
\bordermatrix{
C_{xy}\ & & {\rm column}\ v\in V_2{\setminus}\{y\} &   \cr
{\rm row}\ C_{xy}(u)\ & \ \ldots & \mu_{T{\setminus}xy}(u,v) & \ldots\  \cr }
\end{equation}
By the induction hypothesis applied this time to the tree $T{\setminus}xy$, the row (\ref{eq:13}) is a linear combination of rows $A_{xy}(w)$ of $A_x$, that is, $C_{xy}(u)=\sum_{w\in V_1{\setminus}\{x\}} \varepsilon_w A_{xy}(w)$, with the summation restriction as explained above. Now, each row $A_{xy}(w)$ for $w\in V_1{\setminus} \{x\}$, of length $n-2$, differs from the corresponding row $A(w)$, of length $n$, in the penultimate coordinate (indexed by $y$), of which we do not need to have control, and in the last coordinate (indexed by $x$) equal to $0$. Thus, taking the rows $A(w)$ of the matrix $A$ instead of the rows $A_{xy}(w)$ with the {\em same} coefficients $\varepsilon_w$ the linear combination $S_{xy}(u)=\sum_{w\in V_1{\setminus}\{x\}}\varepsilon_wA(w)$ gives
\begin{equation}\label{eq:14}
\bordermatrix{
 & & {\rm column}\ v\in V_2{\setminus}\{y\} & & {\rm column}\ y \ & \ {\rm column}\ x \  \cr
S_{xy}(u)\ & \ \ldots & \mu_{T{\setminus}xy}(u,v) & \ldots &  \gamma(u,y) & 0 \cr }
\end{equation}
for some integer $\gamma(u,y)$. The last vector we need in this part is the $-\gamma(u,y)$ multiple of the last row $A(x)$ of the matrix $A$ (the row indexed by $x$), which has all coordinates equal to zero except the $(n-1)$-st (indexed by $y$), equal to $1$. It is now obvious from (\ref{eq:12}) and (\ref{eq:14}) that $\theta{\cdot}S_x(u) + S_{xy}(u)- \gamma(u,y)A(x)$ is equal to the $u$-th row $C(u)$ of $C$ displayed in (\ref{eq:10}). Moreover, we have established by induction that $C(u)= \theta{\cdot}S_x(u) + S_{xy}(u)- \gamma(u,y)A(x)$ is a linear combination of rows of $A$, and we proved this for every $u\in V_1$ such that $u\ne x$. If $u=x$, recall the vertex $z$ adjacent to $y$ and consider the last row $C_{xy}(z)$, indexed by $z$, of the matrix $C_{xy}$, which is
\begin{equation}\label{eq:15}
\bordermatrix{
C_{xy} & & {\rm column}\ v\in V_2{\setminus}\{y\} &  \cr
{\rm row}\ C_{xy}(z)\ & \ \ldots & \mu_{T{\setminus}xy}(z,v) & \ldots \ \cr }
\end{equation}
By our induction hypothesis applied $T{\setminus}x$, the row (\ref{eq:15}) can be expressed as a linear combination of rows $A_{xy}(w)$ of $A_x$ in the form $C_{xy}(z)=\sum_{w\in V_1{\setminus}\{x\}} \nu_w A_{xy}(w)$. As in the previous case, observe that each row $A_{xy}(w)$ for $w\in V_1{\setminus} \{x\}$, of length $n-2$, differs from the corresponding row $A(w)$, of length $n$, in the penultimate $y$-th coordinate with value $0$ or $1$ and in the last $x$-th coordinate equal to $0$. Replacing rows $A_{xy}(w)$ with $A(w)$ in the above linear combination but keeping its coefficients $\nu_w$ we obtain $S_{xy}(z)=\sum_{w\in V_1{\setminus}\{x\}}\nu_wA(w)$ where
\begin{equation}\label{eq:16}
\bordermatrix{
 & & {\rm column}\ v\in V_2{\setminus}\{y\} & & {\rm column}\ y \ & \ {\rm column}\ x \ \cr
S_{xy}(z)\ & \ \ldots & \mu_{T{\setminus}xy}(z,v) & \ldots &  \omega(u,y) & 0 \cr }
\end{equation}
for some integer $\omega(u,y)$. And, as in the previous case, the row $C(x)$ of (\ref{eq:10}) is obtained from (\ref{eq:16}) as $C(x)=-S_{xy}(z) + (\omega(u,y)+\mu_T(x,y))A(x)$, which is (by induction) a linear combination of rows of $A$. This completes the analysis in the case when $u\in V_1$. (We reiterate that in all the above diagrams representing rows of matrices, entries in the columns indexed by elements of $V_1{\setminus}\{x\}$ are equal to zero.)
\smallskip

It remains to consider the case when $u\in V_2$; recall that  $y\in V_2$ as well and $x,z\in V_1$. Repeating the previous analysis almost verbatim but interchanging the roles of $V_2$ and $V_1$, that is, considering rows of $C$ indexed by elements $u\in V_2$ and columns indexed by $v\in V_1$ one finds that the rows $C(u)$ for $u\in V_2{\setminus}\{y\}$ and $C(y)$ have the form
\begin{equation}\label{eq:17}
\bordermatrix{
C & & {\rm column}\ v\in V_1{\setminus}\{x\} & & {\rm column}\ y \ & \ {\rm column}\ x \ \cr
{\rm row}\ C(u)\ & \ \ldots & \theta\mu_{T{\setminus}x}(u,v)+\mu_{T{\setminus}xy}(u,v) & \ldots & 0 & -\mu_{T{\setminus}xy}(u,z)   \cr
{\rm row}\ C(y)\ & \ \ldots & \theta\mu_{T{\setminus}x}(y,v) & \ldots & 0 & \mu_T(y,x)\ \cr }
\end{equation}
where we again assume zeros in columns indexed by elements of $V_2$. It should now be clear that one may use induction almost exactly as in the previous case (considering rows of the matrices $C_x$ and $C_{xy}$ indexed by elements of $V_2$) to argue that the rows of $C$ displayed in (\ref{eq:17}) are linear combinations of rows of $A$; we leave the obvious details to the reader. This completes the proof of the statement that the matrix $C$ is Gaussian equivalent to the adjacency matrix $A$ of $T$.
\ebox

\section{Conclusion}\label{sec:fin}

It remains to put the pieces together to prove our main result.
\medskip

{\bf Proof of Theorem \ref{main}.} Let $f$ be an eigenvector of $T$ corresponding to a non-zero eigenvalue. By Propositions \ref{prop:satur} and \ref{prop:unsat}, the same $f$ satisfies the part (a) of Proposition \ref{prop:eigen}. Further, by Proposition \ref{prop:Gauss} the matrices $A$ and $B=m(T)^{-1}C$ are Gaussian equivalent and hence have the same null-spaces. It follows that the matrix $A$ with all its eigenvectors, and the matrix $B$ indexed the same way as $A$ with entries $\mu_T(u,v)/m(T)$, satisfy the assumptions of Proposition \ref{prop:eigen}, by which $B$ is the generalized inverse of $A$. Thus, $B$ is an adjacency matrix of the inverse tree $(T^*,\alpha^*)$, completing the proof. \ebox

Our approach based on Proposition \ref{prop:eigen} requires having a preliminary idea about entries of a generalized inverse. In our case, wanting to extend the inversion formula of \cite{PK} for trees with a (unique) perfect matching, a further educated guess lends itself by considering trees with simple structure but containing a large number of maximum matchings, such as paths of even lengths or stars. Taking on the latter, if $T$ is a star with $n$ pendant edges, its adjacency matrix $A$ has spectrum $(\pm\sqrt{n},0^{n-1})$ and it is easy to check that its generalized inverse is $A^*=n^{-1}A$, with spectrum $(\pm 1/\sqrt{n},0^{n-1})$. This immediately yields the clue of including the term $m(T)=n$ in the denominator of the formula of \cite{PK}, leading eventually to Theorem \ref{main}. 
\smallskip

We believe that our new method of considering {\em eigenspaces in an implicit form} will lead to extensions of Theorem \ref{main} to broader classes of graphs, at the very least in a similar way Godsil's inversion theorem for trees with a perfect matching \cite{Gods} and the subsequent formula for inverses of such trees \cite{PK} have been extended to bipartite graphs with a unique perfect matching in \cite{BNP,AK,Bapat,McMc,Pavl}.

\bigskip

\noindent{\bf Acknowledgment}\ \ Both authors acknowledge support by research grants APVV 0136/12 and APVV-15-0220, and by research grants VEGA 1/0026/16 and VEGA 1/0142/17.

\end{document}